\documentclass[11pt]{amsart}
\pdfoutput=1
\usepackage{amssymb,amscd,graphicx,enumerate,epsfig}
\usepackage[bookmarks=true]{hyperref}

\newcommand{\Claim}[1]{\noindent\textbf{Claim. }\textit{#1}\\}
\newtheorem{theorem}{Theorem}[section]

\newtheorem{proposition}[theorem]{Proposition}

\newtheorem{lemma}[theorem]{Lemma}

\theoremstyle{definition}

\newenvironment{proofc}{\noindent\textit{Proof of Claim.}}{\\}

\newcommand{\Case}[1]{\textbf{Case #1.}}
\begin{document}

\title[Finite group actions on Haken Manifolds]{Handlebody-preserving finite group actions on Haken manifolds with Heegaard genus two - II}
\author{Jungsoo Kim}
\date{}
\begin{abstract} Let $M$ be a closed orientable $3$-manifold of Heegaard genus two with a non-trivial JSJ-decomposition and $G$ be a finite group of orientation-preserving diffeomorphisms acting on $M$ which preserves each handlebody of Heegaard splitting and each piece of the JSJ-decomposition of $M$. Let $(V_1,V_2;F)$ be the Heegaard splitting and $\cup T_i$ be the union of the JSJ-tori. In this article, we prove that $G \cong \mathbb{Z}_2$ or $\mathbb{D}_2$ if all components of $V_j\cap (\cup T_i)$ are not $\partial$-parallel in $V_j$ for $j=1,2$ unless $V_j\cap(\cup T_i)$ has three or more disk components for $j=1$ or $2$.
\end{abstract}

\address{\parbox{4in}{Department of Mathematics\\
Konkuk University\\
Seoul 143-701\\ Korea\medskip}} \email{pibonazi@gmail.com}

\maketitle
\tableofcontents

\section[Introduction]{Introduction}
Let $M$ be a closed orientable $3$-manifold of Heegaard genus two with a non-trivial JSJ-decomposition and $G$ be a finite group of orientation-preserving diffeomorphisms acting on $M$ which preserves each handlebody of Heegaard splitting and each piece of the JSJ-decomposition of $M$. Let $(V_1,V_2;F)$ be the Heegaard splitting and $\cup T_i$ be the union of the JSJ-tori.

In \cite{KJS2}, there are two conditions \textit{``Condition A''} and \textit{``Condition B''} for the Heegaard splitting and the JSJ-tori, where each determines $G$ to be $\mathbb{Z}_2$ or $\mathbb{D}_2$ up to isomorphism. But these two conditions have strong restrictions - $V_1\cap (\cup T_i)$ must consist of at most two disks or annuli. How can we generalize these conditions?

In this article, we will determine the possible isomorphism types of $G$ by the following theorem.

\begin{theorem}[the Main theorem]\label{main-theorem}
Let $M$ be a closed orientable $3$-manifold with a genus two Heegaard splitting $(V_1, V_2; F)$ and a non-trivial JSJ-decomposition, where all components of the intersection of the JSJ-tori and $V_i$ are not $\partial$-parallel in $V_i$ for $i=1,2$. If $G$ is a finite group of orientation-preserving smooth actions on $M$ which preserves each handlebody of the Heegaard splitting and each piece of the JSJ-decomposition of $M$, then $G\cong \mathbb{Z}_2$ or $\mathbb{D}_2$ unless $V_j\cap(\cup T_i)$ has three or more disk components for $j=1$ or $2$.
\end{theorem}

\begin{center}
Acknowledgement
\end{center}
This work was supported by the Korea Science and Engineering
Foundation (KOSEF) grant funded by the Korea government (MOST)
(No.~R01-2007-000-20293-0).

\section{Proof of the Main Theorem}

We define two conditions for the Heegaard splitting and the JSJ-tori as follows.
\begin{enumerate}[(A).]
    \item $V_1\cap(\cup T_i)$ consists of disks, and the number of disks is at most two,
    \item $V_1\cap(\cup T_i)$ consists of annuli, and the number of annuli is at most two,\\
where we assume that each component of $V_1\cap(\cup T_i)$ is not $\partial$-parallel in $V_1$.
\end{enumerate}
Let us call the first \textit{Condition A} and the second \textit{Condition B} (see section 1 of \cite{KJS2} for more details.)

In \cite{KJS2}, we get the result under Condition A or Condition B.

\begin{proposition}\label{theorem-KJS2}(\cite{KJS2}, Theorem 1.1)
Let $M$ be a closed orientable $3$-manifold with a genus two Heegaard splitting $(V_1, V_2; F)$ and a non-trivial JSJ-decomposition, where all components of the intersection of the JSJ-tori and $V_i$ are not $\partial$-parallel in $V_i$ for $i=1,2$. If $G$ is a finite group of orientation-preserving diffeomorphisms acting on $M$ which preserves each handlebody of the Heegaard splitting and each piece of the JSJ-decomposition of $M$, then $G\cong \mathbb{Z}_2$ or $\mathbb{D}_2$ if $V_j\cap(\cup T_i)$ consists of at most two disks or at most two annuli.
\end{proposition}

So now we assume that the Heegaard splitting and the JSJ-tori satisfy neither Condition A nor Condition B.

Here we introduce some geometric lemmas for the proof of the Main theorem.

\begin{lemma}(\cite{KJS2}, Lemma 2.8)\label{lemma-gm-12}
Let  $G$ be a finite group of non-trivial orientation preserving
diffeomorphisms which acts on a genus two handlebody $V$. Suppose that $\mathcal{S}$ is a set of properly embedded surfaces in $V$ where $\mathcal{S}$ cuts $V$ into two solid tori $V_1$, $V_3$ and a $3$-ball $V_2$. If $G$ preserves each $V_i$ for $i=1,2,3$, or $G$ preserves $V_2$ but exchanges $V_1$ and $V_3$, then $G\cong\mathbb{Z}_2$ or $\mathbb{D}_2$.
\end{lemma}

\begin{lemma}(\cite{KJS2}, Lemma 2.9)\label{lemma-gm-2}
Let  $G$ be a finite group of non-trivial orientation preserving diffeomorphisms which acts on a genus two handlebody $V$. Suppose that $\mathcal{S}$ is a set of properly embedded surfaces in $V$ where $\mathcal{S}$ cuts $V$ into a $3$-ball $V_1$ and a solid torus $V_2$. If $G$ preserves both $V_1$ and $V_2$, then $G\cong  \mathbb{Z}_2$ or $\mathbb{D}_2$.
\end{lemma}

\begin{lemma}(\cite{KJS2}, Lemma 2.11)\label{lemma-gm-an-1}
Let $G$ be a finite group of non-trivial orientation preserving diffeomorphisms which acts on a genus two handlebody $V$. Suppose that $\mathcal{S}$ is a set of properly embedded surfaces in $V$ where $\mathcal{S}$ cuts $V$ into two connected components $V_1$ and $V_2$, where both $V_1$ and $V_2$ are noncontractible in $V$, and $V_1$ or $V_2$ contains a meridian disk of $V$. If $G$ preserves both $V_1$ and $V_2$, then $G\cong  \mathbb{Z}_2$ or $\mathbb{D}_2$.
\end{lemma}

\begin{lemma}(\cite{KJS2}, Lemma 2.12)\label{lemma-gm-an-2}
Let  $G$ be a finite group of non-trivial orientation preserving diffeomorphisms which acts on a genus two handlebody $V$. Suppose that $\mathcal{S}$ is a set of properly embedded surfaces in $V$ where $\mathcal{S}$ cuts $V$ into two solid tori $V_1$ and $V_3$, and a genus two handlebody $V_2$, where
    \begin{enumerate}
        \item $V'=\overline{V_1}\cup \overline{V_2}$ is a genus two handlebody,
        \item both $V_1$ and $V_2$ are non-contractible in $V'$ and $V$, and
        \item there is a meridian disk $D$ in $V$ such that $D\cap V_1=\emptyset$ and $D\cap V_2$ is a non-empty set of meridian disks of $V'$.
    \end{enumerate}
If $G$ preserves $V_i$ for $i=1,2,3$, then $G\cong  \mathbb{Z}_2$ or $\mathbb{D}_2$. If $G$ preserves $V_2$ and exchanges $V_1$ and $V_3$, then $G\cong \mathbb{Z}_2$ or $\mathbb{D}_2$.
\end{lemma}

\begin{figure}
\centering
\includegraphics[viewport=31 128 508 780, width=8cm]{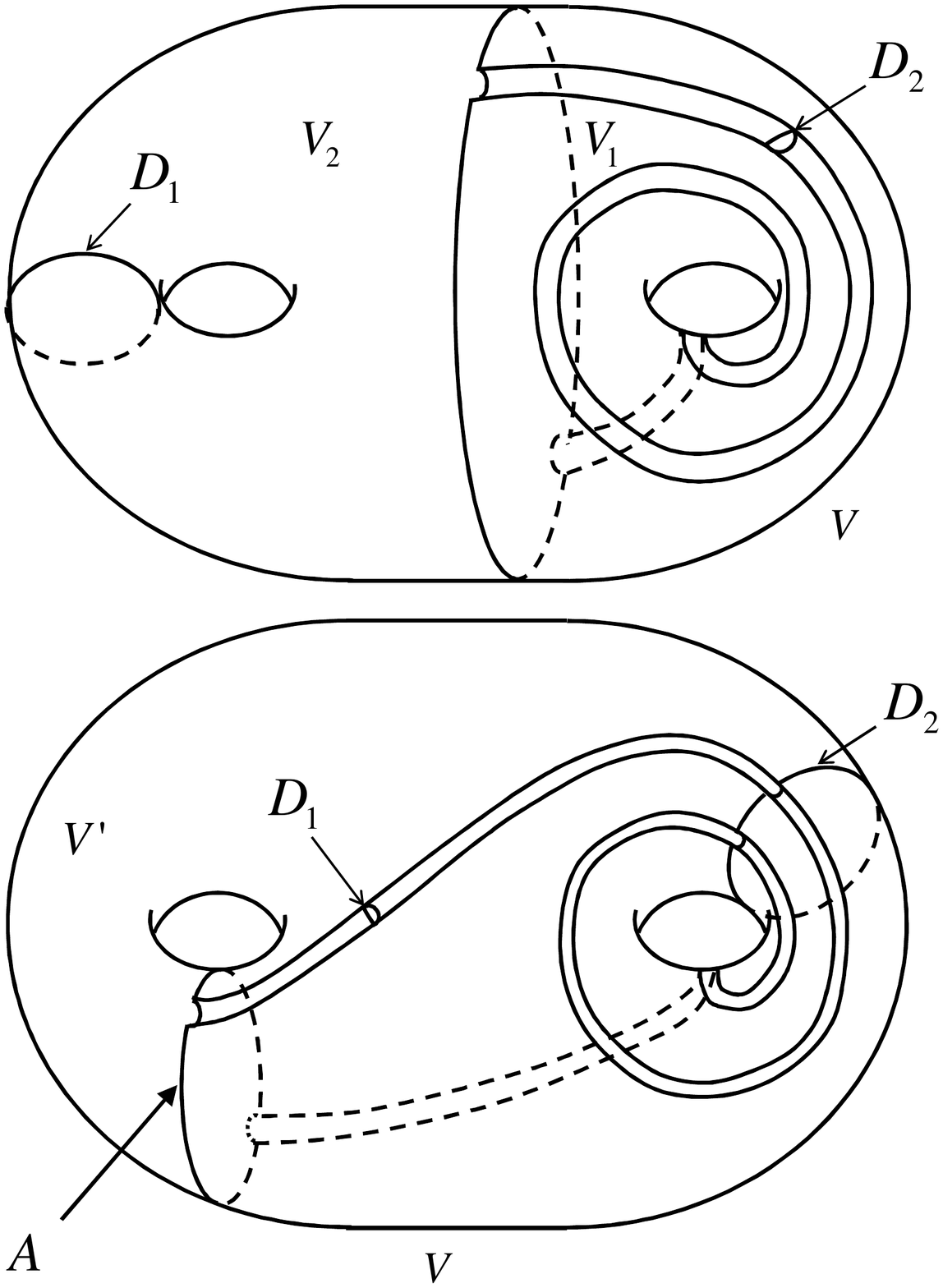}
\caption{Two possible positions of an essential annulus in a genus two handlebody. \cite{KO}\label{fig-lemma-ko-3-2}}
\end{figure}

\begin{lemma}\label{lemma-ko-3-2}(\cite{KO}, Lemma 3.2) If $A$ is an essential annulus in a genus two handlebody $V$ then either
    \begin{enumerate}[(i)]
        \item $A$ cuts $V$ into a solid torus $V_1$ and a genus two handlebody $V_2$ and there is a complete system of meridian disks $\{D_1,D_2\}$ of $V_2$ such that $D_1\cap A=\emptyset$ and $D_2\cap A$ is an essential arc of $A$ (see the upper of Figure \ref{fig-lemma-ko-3-2},) or\label{lemma-ko-3-2-1}
        \item $A$ cuts $V$ into a genus two handlebody $V'$ and there is a complete system of meridian disks $\{D_1,D_2\}$ of $V'$ such that $D_1\cap A$ is an essential arc of $A$ (see the lower of Figure \ref{fig-lemma-ko-3-2}.)\label{lemma-ko-3-2-2}
    \end{enumerate}

\end{lemma}

\begin{figure}
\centering
\includegraphics[viewport=26 231 416 780, width=7cm]{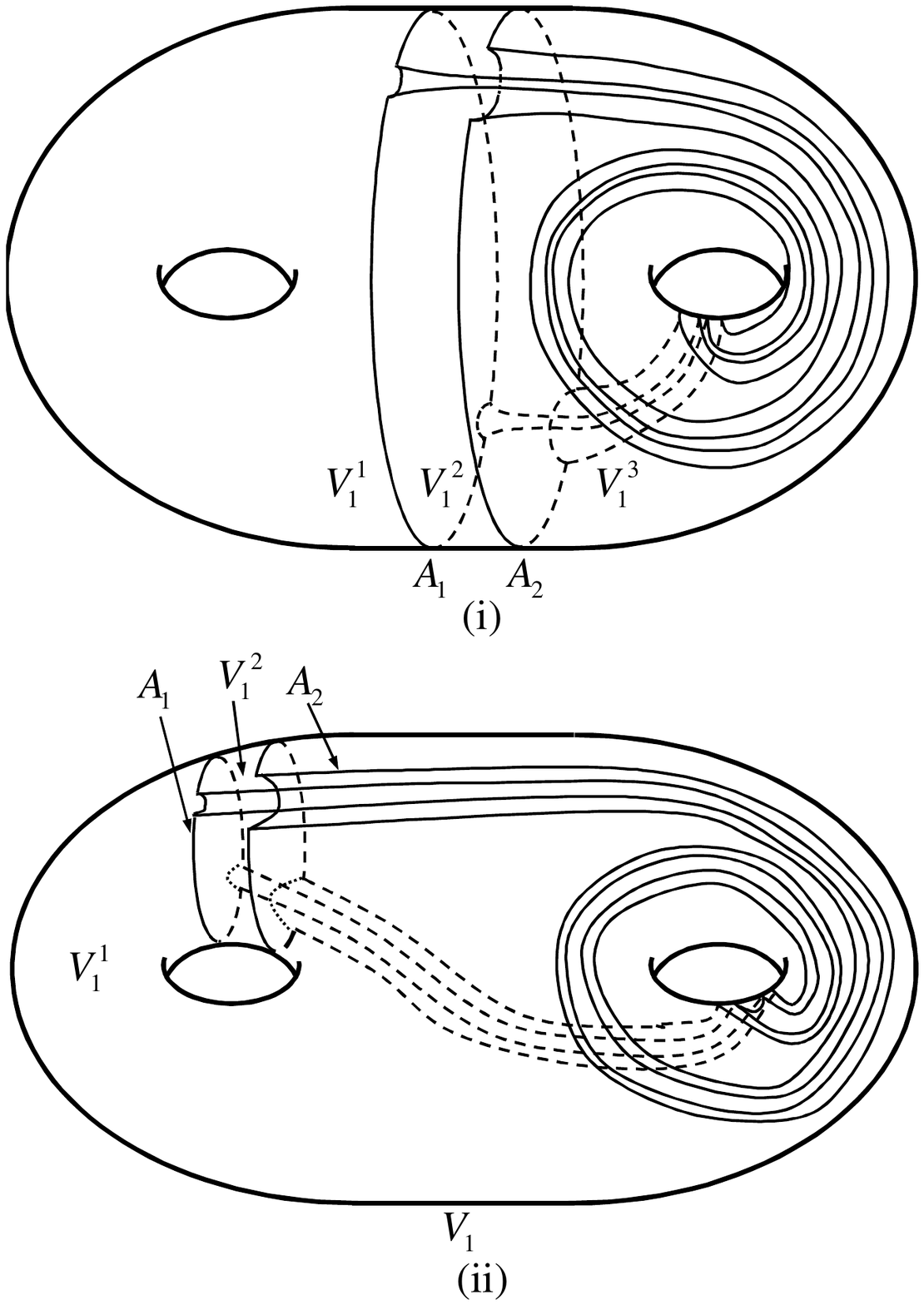}
\caption{$A_1$ and $A_2$ are parallel in $V_1$.\label{figure-add-ad3}}
\end{figure}

In particular, if two essential annuli $A_1$ and $A_2$ are parallel in a genus two handlebody $V_1$, then we get Figure \ref{figure-add-ad3}.

\begin{figure}
\centering
\includegraphics[viewport=28 144 336 777, width=8cm]{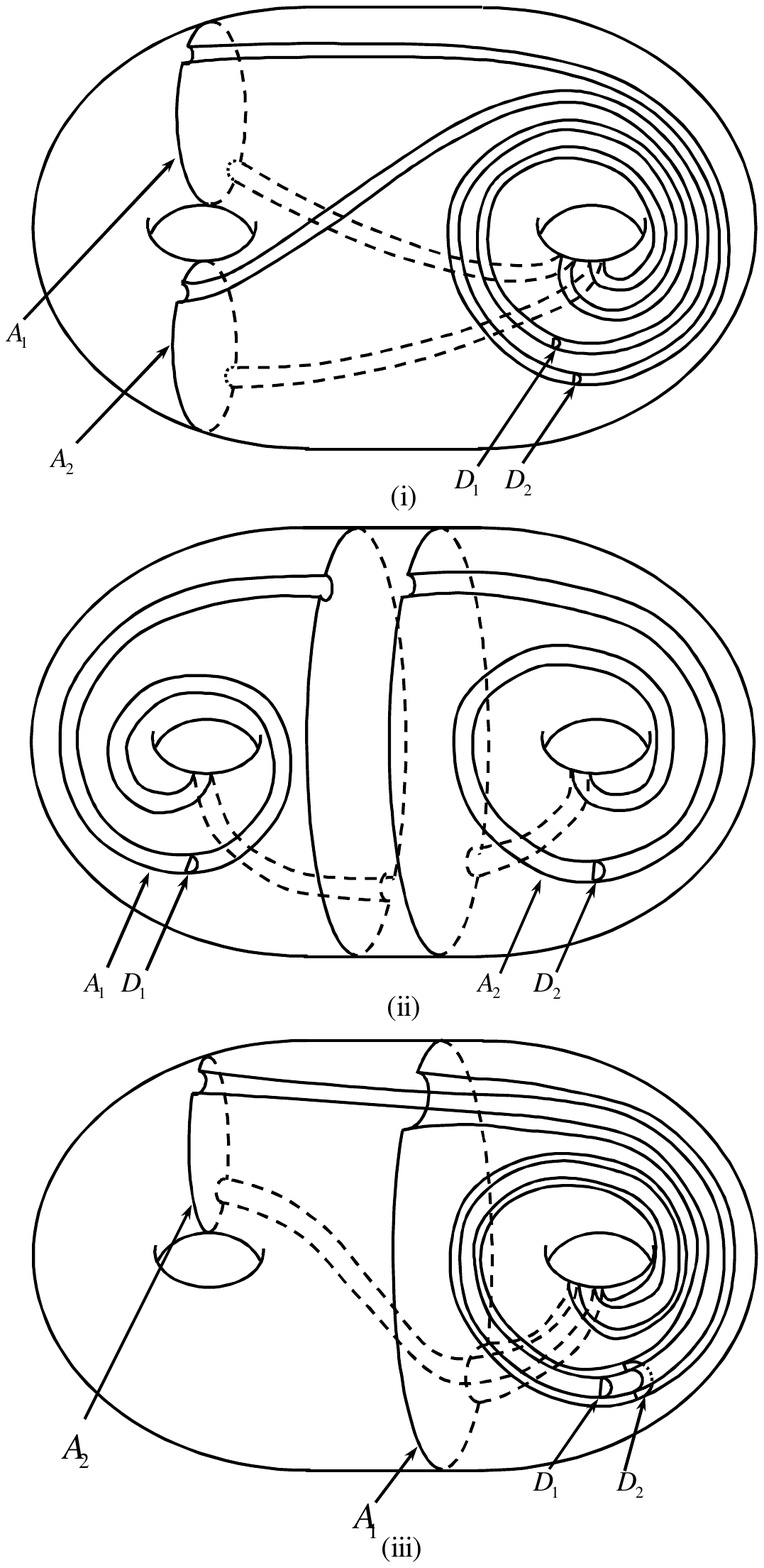}
\caption{The possibilities of two mutually disjoint non-parallel essential annuli in the genus two handlebody. \cite{KO}\label{fig-lemma-ko-3-4}}
\end{figure}

\begin{lemma}\label{lemma-ko-3-4}(\cite{KO}, Lemma 3.4) Let $\{A_1, A_2\}$ be a system of mutually disjoint, non-parallel, essential annuli in a genus two handlebody $V$. Then either
    \begin{enumerate}[(i)]
        \item \label{lemma-ko-3-4-1}$A_1\cup A_2$ cuts $V$ into a  solid torus $V_1$ and a genus two handlebody $V_2$. Then $A_1\cup A_2\subset \partial V_1$, $A_1\cup A_2\subset\partial V_2$ and there is a complete system of meridian disks $\{D_1,D_2\}$ of $V_2$ such that $D_i\cap A_j=\emptyset$ ($i\neq j$) and $D_i\cap A_i$ ($i=1,2$) is an essential arc of $A_i$,
        \item \label{lemma-ko-3-4-2}$A_1\cup A_2$ cuts $V$ into two solid tori $V_1$, $V_2$ and a genus two handlebody $V_3$. Then $A_1\subset \partial V_1$, $A_2\subset \partial V_2$, $A_1\cup A_2 \subset \partial V_3$ and there is a complete system of meridian disks $\{D_1, D_2\}$ of $V_3$ such that $D_i\cap A_j=\emptyset$ $(i\neq j)$ and $D_i\cap A_i$ ($i=1,2$) is an essential arc of $A_i$ or
        \item \label{lemma-ko-3-4-3}$A_1\cup A_2$ cuts $V$ into a solid torus $V_1$ and a genus two handlebody $V_2$. Then $A_i\subset\partial V_1$ ($i=1$ or $2$, say $1$), $A_2\cap V_1=\emptyset$, $A_1\subset\partial V_2$ and there is a complete system of meridian disks $\{D_1, D_2\}$ of $V_2$ such that $D_1\cap A_2$ is an essential arc of $A_2$ and $D_2\cap            A_i$ ($i=1,2$) is an essential arc of $A_i$, see Figure \ref{fig-lemma-ko-3-4}.
    \end{enumerate}
\end{lemma}

Since a JSJ-torus is incompressible in $M$, all components of $V_1\cap (\cup T_i)$ are incompressible in $V_1$. So we will assume the incompressibility of all components of $V_1\cap (\cup T_i)$ in $V_1$. By the assumption of Theorem \ref{main-theorem}, we assume that all components of $V_1\cap (\cup T_i)$ are not $\partial$-parallel in $V_1$.

At first, we will prove that $G\cong \mathbb{Z}_2$ or $\mathbb{D}_2$ if $V_1\cap (\cup T_i)$ contains at most two disk components by Lemma \ref{lemma-d1cp} and Lemma \ref{lemma-d2cp3morecp}.

\begin{lemma}\label{lemma-d1cp}
If $V_1\cap (\cup T_i)$ has only one disk component $D$, then $G\cong \mathbb{Z}_2$ or $\mathbb{D}_2$.
\end{lemma}

\begin{proof}
Since $G$ preserves $V_1\cap (\cup T_i)$ and there is only one disk component $D$ in $V_1\cap (\cup T_i)$, $G$ must preserve $D$. So $G$ also preserves a small product neighborhood of $D$ in $V_1$, say $\overline{N(D)}\cong D^2 \times I$, where $\overline{N(D)}\cap \partial V_1=\partial D^2\times I$ and $N(D)$ meets no component of $V_1\cap (\cup T_i)$ other than $D$. Say that $D_1=D^2\times\{0\}$ and $D_2=D^2\times\{1\}$.
Then $\{D_1, D_2\}$ cuts $V_1$ into a $3$-ball $B$ and a solid torus (when $D_1$ is non-separating in $V_1$,) or into a $3$-ball $B$ and two solid tori (when $D_1$ is separating in $V_1$.) Since $G$ preserves $\{D_1,D_2\}$ and  $B$ is the only $3$-ball component of $V_1-\{D_1,D_2\}$, $G$ preserves $B$. So $G\cong \mathbb{Z}_2$ or $\mathbb{D}_2$ by Lemma \ref{lemma-gm-2} or Lemma \ref{lemma-gm-12}.
\end{proof}

\begin{lemma}\label{lemma-d2cp3morecp}
If $V_1\cap (\cup T_i)$ consists of at least three components and has exactly two disk components $D_1$ and $D_2$, then $G\cong \mathbb{Z}_2$ or $\mathbb{D}_2$.
\end{lemma}

\begin{figure}
\centering
\includegraphics[viewport=37 400 378 778, width=6cm]{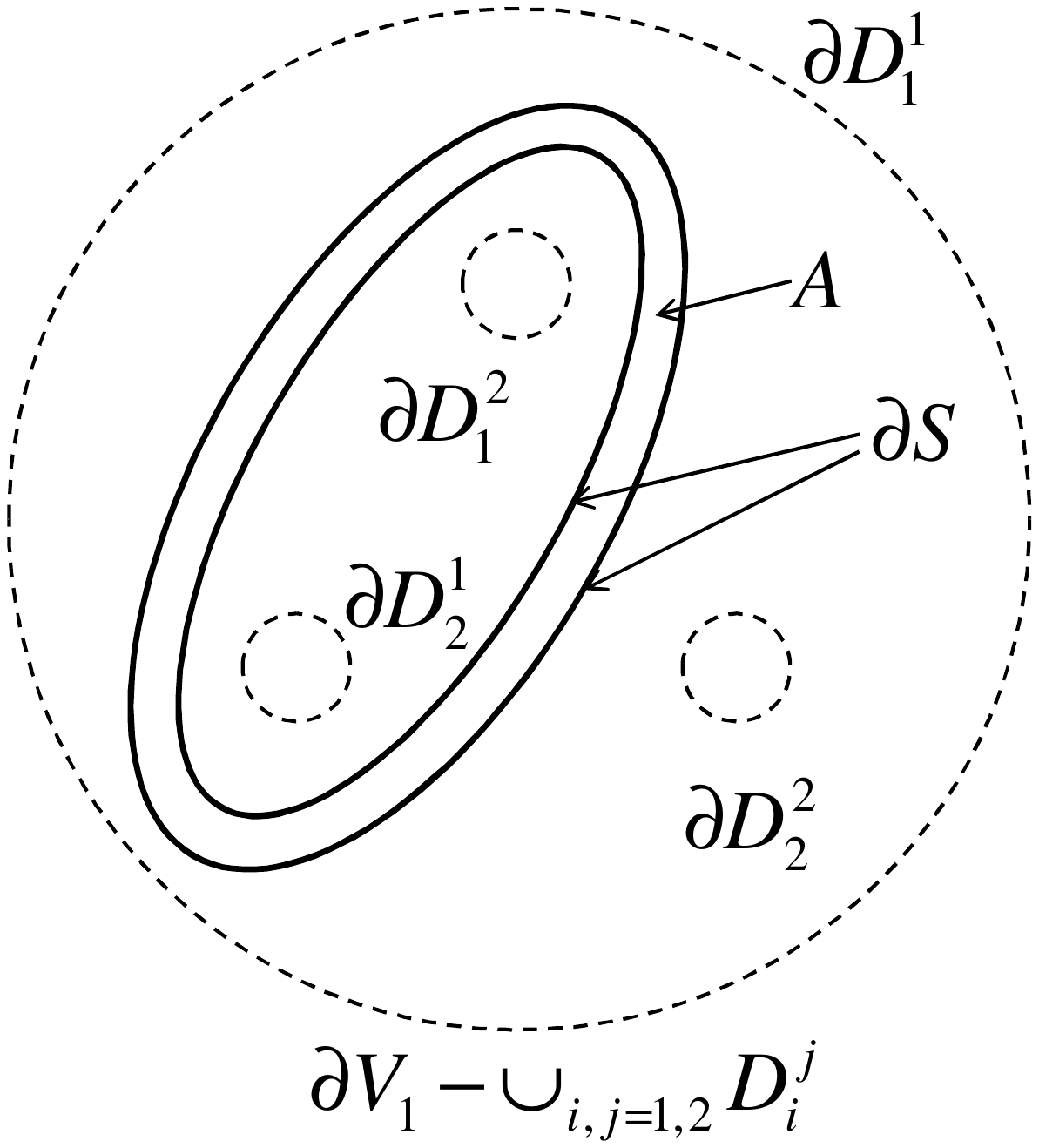}
\caption{$S$ is $\partial$-parallel to an annulus in $\partial V_1$.\label{Bprime}}
\end{figure}

\begin{proof}
We divide the proof into two cases.\\

\Case{1} $D_1$ and $D_2$ are parallel in $V_1$.

Since $D_1$ and $D_2$ are the only disk components in $V_1\cap (\cup T_i)$, $G$ must preserve $\{D_1, D_2\}$ in $V_1$.
So $G$ preserves $V_1-\{D_1,D_2\}$. Since the $3$-ball $B$ between $D_1$ and $D_2$ is the only $3$-ball component of $V_1-\{D_1,D_2\}$ (see $D_1$ and $D_2$ in the proof of Lemma \ref{lemma-d1cp}), $G$ preserves $B$. So $G\cong \mathbb{Z}_2$ or $\mathbb{D}_2$ using the arguments in the proof of Lemma \ref{lemma-d1cp}.\\

\Case{2} $D_1$ and $D_2$ are non-parallel in $V_1$.

If one of $D_1$ and $D_2$, say $D_1$, is separating in $V_1$, then $D_2$ is non-separating in $V_1$. So $\{D_1, D_2\}$ cuts $V_1$ into a $3$-ball $V_1^1$ and a solid torus $V_1^2$. Since $G$ preserves $\{D_1, D_2\}$, $G$ also preserves $V_1-\{D_1, D_2\}=V_1^1\cup V_1^2$. But $V_1^1$ is not homeomorphic to $V_1^2$, $G$ must preserves both $V_1^1$ and $V_1^2$. So $G\cong \mathbb{Z}_2$ or $\mathbb{D}_2$ by Lemma \ref{lemma-gm-2}.

If both of $D_1$ and $D_2$ are non-separating in $V_1$, then $\{D_1,D_2\}$ cuts $V_1$ into a $3$-ball $B$. Let $D_i^1$ and $D_i^2$ be two copies of $D_i$ for $i=1,2$ in $\partial B$ after cutting $V_1$ along $D_1$ and $D_2$, and $S$ be a non-disk component of $V_1\cap (\cup T_i)$. Since any properly embedded incompressible surface in $B^3$ must be $\partial$-parallel, $S$ is $\partial$-parallel in $B$. Moreover, no curve of $\partial S$ is parallel to some $\partial D_i^j$ for $1\leq i,j\leq 2$ ($S$ is incompressible in $V_1$ and does not intersect both $D_1$ and $D_2$.) Let $B'$ be $\partial B-\cup_{i=1,2} D_i^j$. Then $B'$ is a disk with three holes. So each curve of $\partial S$ must separate each two of $\{\partial D_i^j\}_{i=1,2}$. In addition, the curves of $\partial S$ must be pairwise parallel in $B'$ after one of them is determined (see Figure \ref{Bprime}.) Since $S$ is a non-disk connected surface and all curves in $\partial S$ are pairwise parallel in $B'$, $\partial S$ must consist of two curves. So $S$ is $\partial$-parallel to an annulus $A\subset B'$ in $B$, where $A$ is bounded by these two curves, i.e. $S$ is $\partial$-parallel in $V_1$. But this contradicts the hypothesis.
\end{proof}

Now we assume that $V_1\cap (\cup T_i)$ consists of non-disk components.
We will prove that $G\cong \mathbb{Z}_2$ or $\mathbb{D}_2$ if $V_1\cap (\cup T_i)$ contains a non-annuls component.

\begin{lemma}\label{lemma-nonannuluscp}
Suppose that $V_1\cap (\cup T_i)$ consists of non-disk components.
If $V_1\cap (\cup T_i)$ has a non-annulus component, then $G\cong \mathbb{Z}_2$ or $\mathbb{D}_2$ unless $V_2\cap (\cup T_i)$ has three or more disk components.
\end{lemma}

\begin{proof}
Let $S$ be a non-annulus component of $V_1\cap (\cup T_i)$, $T_1$ be a JSJ-torus which contains $S$.
Then $S$ is a disk with $n$ holes ($n\geq 2$) or a torus with $m$ holes ($m\geq 1$) since $S\subset T_1$.\\

\Claim{$T_1-\partial V_1$ must contain at least one disk component.}

    \begin{proofc}
    $T_1-S$ has at least one disk component $D$. If we choose an innermost disk of $D-\partial V_1$, then the proof ends.
    \end{proofc}

So $V_2\cap (\cup T_i)$ contains at least one disk component by the above claim. Since $G$ also preserves $V_2$, if we substitute $V_2$ for $V_1$, then we get $G\cong \mathbb{Z}_2$ or $\mathbb{D}_2$ by Lemma \ref{lemma-d1cp} or Lemma \ref{lemma-d2cp3morecp} unless $V_2\cap (\cup T_i)$ has three or more disk components.
\end{proof}

Now we remains only the cases when $V_1\cap (\cup T_i)$ consists of annuli. Since we already proved the cases when $V_1\cap (\cup T_i)$ consists of at most two annulus components in ``Condition B'', it is sufficient to prove the cases $V_1\cap (\cup T_i)$ consists of at least three annulus components.

At first, we will prove $V_1\cap (\cup T_i)$ consists of pairwise non-parallel annuli.

T. Kobayashi proved that there is a concrete way to express the system of pairwise disjoint, non-parallel three essential annuli in the genus two handlebody by the following lemma.

\begin{figure}
\centering
\includegraphics[viewport=29 446 525 780, width=8cm]{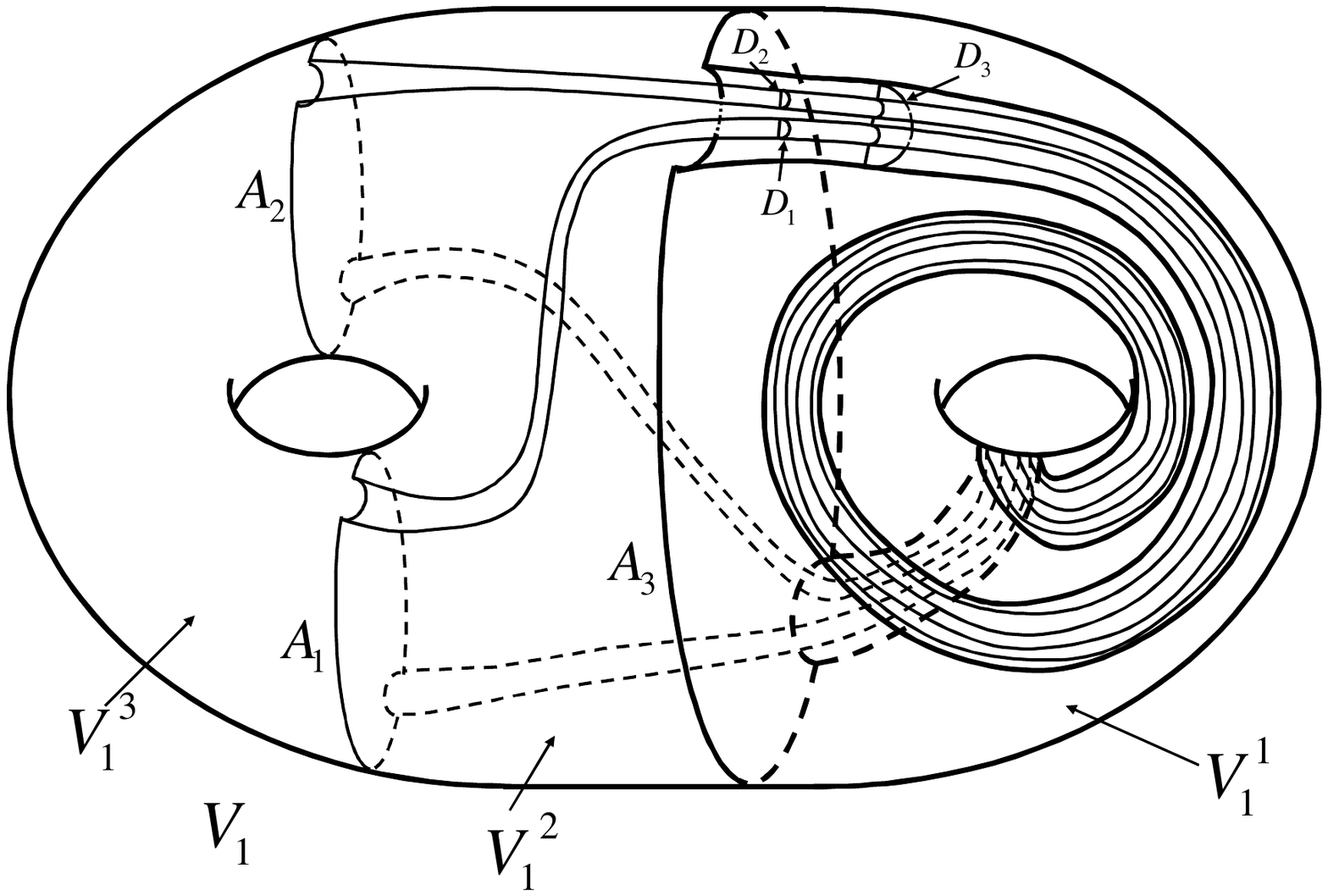}
\caption{$\{A_1,A_2,A_3\}$ cuts $V_1$ into a genus two handlebody and two solid tori.\label{fig-threeannuli}}
\end{figure}

\begin{lemma}(\cite{KO}, Lemma 3.5)\label{lemma-ko-3-5}
Let $\{A_1,A_2,A_3\}$ be a system of pairwise disjoint, non-parallel essential annuli in the genus two handlebody $V_1$. Then $A_1\cup A_2 \cup A_3$ cuts $V_1$ into two solid tori $V_1^1$, $V_1^2$ and a genus two handlebody $V_1^3$ which satisfies
    \begin{enumerate}
            \item $A_i\subset \partial V_1^1$ ($i=1$, $2$ or $3$, say $3$), $A_1$, $A_2\subset \partial V_1^3$, $A_1$, $A_2$, $A_3\subset \partial V_1^2$,
            \item there is a complete system of meridian disks $\{D_1, D_2\}$ of $V_1^3$ such that $D_i\cap A_j=\emptyset$ $(i\neq j)$ and $D_i\cap A_i$ ($i=1,2$) is an essential arc of $A_i$ and
            \item there is a meridian disk $D_3$ of $V_1^2$ such that $D_3\cap A_i$ ($i=1,2,3$) is an essential arc of $A_i$.
    \end{enumerate}
See Figure \ref{fig-threeannuli}.
\end{lemma}

\begin{lemma}\label{lemma-3moreannulicp}
If $V_1\cap (\cup T_i)$ consists of three or more pairwise non-parallel annulus components, then $G\cong \mathbb{Z}_2$ or $\mathbb{D}_2$.
\end{lemma}

\begin{proof}

At first, we introduce the following claim.

\Claim{
There is no system of properly embedded, pairwise disjoint, non-parallel essential four or more annuli in the genus two handlebody $V_1$.}

\begin{proofc}
Suppose that there is a system of pairwise disjoint, non-parallel essential four annuli in a genus two handlebody, say $\{A_1, A_2, A_3, A_4\}$. $\{A_1, A_2, A_3\}$ cuts $V_1$ into a genus two handlebody $V_1^3$ and two solid tori $V_1^1$ and $V_1^2$ as Lemma \ref{lemma-ko-3-5}. We know that any properly embedded incompressible annulus in a solid torus must be $\partial$-parallel in the solid torus from Lemma 3.1 of \cite{KO}. If $A_4$ is embedded in $V_1^1$ or $V_1^2$, then $A_4$ must be $\partial$-parallel in $V_1^1$ or $V_1^2$. But $A_4$ cannot be parallel into any one of $A_1$, $A_2$ and $A_3$, $A_4$ must be $\partial$-parallel to a subset of $\partial V_1^1-\operatorname{int}(A_3)$ or $\partial V_1^2-\operatorname{int}(A_1\cup A_2\cup A_3)$ ($\partial A_4$ is in $\partial V_1^1$ or $\partial V_1^2$.) So $A_4$ is $\partial$-parallel in $V_1$, this is a contradiction. So $A_4$ is embedded in $V_1^3$. Moreover, $A_4$ must be satisfy (\ref{lemma-ko-3-2-1}) of Lemma \ref{lemma-ko-3-2} in $V_1^3$ ($A_4$ must be properly embedded not only in $V_1^3$, but also in $V_1$, see Figure \ref{fig-lemma-ko-3-2} and Figure \ref{fig-threeannuli}.)  Say that $\{A_1, A_2\}$ bounds $V_1^3$. Then $A_4$ must be parallel to $A_1$ or $A_2$ in $V_1^3$. ($A_4$ is properly embedded not only in $V_1^3$, but also in $V_1$, too. See Figure \ref{fig-threeannuli}.) So we get a contradiction.
\end{proofc}

So we assume that $V_1\cap (\cup T_i)$ consists of exactly three, non-parallel annuli $A_1$, $A_2$ and $A_3$. Use the same notations of $A_1$, $A_2$, $A_3$, $V_1^1$, $V_1^2$ and $V_1^3$ as Lemma \ref{lemma-ko-3-5}. Since $V_1^3$ is the only genus two handlebody among $V_1^1$, $V_1^2$ and $V_1^3$, $G$ preserves $V_1^3$, i.e. $G$ preserves $\{A_1,A_2\}$. So $G$ also preserves $V_1^2$ (so also preserves $V_1^1$.) Let $V'$ be  $\overline{V_1^3}\cup \overline{V_1^2}$. Then $V'$ is a genus two handlebody and $V_1=V'\cup_{A_3}\overline{ V_1^1}$. Since $G$ preserves $V'$ and $V_1^1$, we get $G\cong \mathbb{Z}_2$ or $\mathbb{D}_2$ by Lemma \ref{lemma-gm-an-1}.
\end{proof}

Now we only need to prove the following lemma.

\begin{lemma}\label{lemma-parallelannulicp}
If $V_1\cap (\cup T_i)$ consists of at least three annulus components and some annuli $A_1$ and $A_2$ are parallel in $V_1$, then $G\cong \mathbb{Z}_2$ or $\mathbb{D}_2$.
\end{lemma}

\begin{proof}
By Lemma \ref{lemma-ko-3-2}, $\{A_1,A_2\}$ cuts $V_1$ into a solid torus and a genus two handlebody (if $A_1$ is non-separating in $V_1$) or two solid tori and a genus two handlebody (if $A_1$ is separating in $V_1$ ) as Figure \ref{figure-add-ad3}.

If some component of $V_1\cap (\cup T_i)$ except $A_1$ and $A_2$ intersects the solid torus bounded by $A_1$ and $A_2$, then it must be parallel to $A_1$ (also to $A_2$) in $V_1$. So if we choose nearest parallel annuli for $A_1$ and $A_2$, then we can assume that no component of $V_1\cap (\cup T_i)$ except $A_1$ and $A_2$ intersects the solid torus bounded by $A_1$ and $A_2$. If $G$ preserves the solid torus between $A_1$ and $A_2$, then we get $G\cong \mathbb{Z}_2$ or $\mathbb{D}_2$ using the same arguments in Case 3 of section 6 in \cite{KJS2}. So now we assume that $G$ does not preserve the solid torus between $A_1$ and $A_2$. Since there is no system of properly embedded, pairwise disjoint, non-parallel essential four or more annuli in a genus two handlebody from Claim of Lemma \ref{lemma-3moreannulicp}, we can divide the proof into three cases as follows.

\Case{1} $V_1\cap (\cup T_i)$ consists of annuli parallel to $A_1$. If the number of components of $V_1\cap (\cup T_i)$ is $n$, then $V_1\cap (\cup T_i)$ cuts $V_1$ into a genus two handlebody $V_1^1$ and $n-1$ solid tori (if $A_1$ is non-separating, see (ii) of Lemma \ref{lemma-ko-3-2}.) or into a genus two handlebody $V_1^1$ and $n$ solid tori (if $A_1$ is separating, see (i) of Lemma \ref{lemma-ko-3-2}.) In any case, the only genus two handlebody component $V_1^1$ of $V_1-(\cup T_i)$ is preserved by $G$. In the former case, the big solid torus $V_1-V_1^1$ must be preserved by $G$, so we get $G\cong \mathbb{Z}_2$ or $\mathbb{D}_2$ by Lemma \ref{lemma-gm-an-1}. In the latter case, each of adjacently attached $n-1$ solid tori among $n$ solid tori meets two annuli of $V_1\cap (\cup T_i)$ in its closure and the last solid torus $V_1^3$ meets only one annulus of $V_1\cap (\cup T_i)$ in its closure (the naming of $V_1^3$ is the same as (i) of Figure \ref{figure-add-ad3}.) Therefore, $V_1^3$ must be preserved by $G$, i.e. the big solid torus $\tilde{V}=V_1-(V_1^1\cup V_1^3)$ must be preserved by $G$. Since $V'=\overline{V_1^1}\cup \tilde{V}$ is a genus two handlebody, we get $G\cong \mathbb{Z}_2$ or $\mathbb{D}_2$ by Lemma \ref{lemma-gm-an-1} (in (i) of Figure\ref{figure-add-ad3}, we can say that the big solid torus $\tilde{V}$ is $V_1^2$.)

\Case{2} The parallel class of $V_1\cap (\cup T_i)$ in $V_1$ consists of two elements.
Let the elements of the class be $[A']$ and $[A'']$. Then $\{A',A''\}$ satisfies one of three conclusions of Lemma \ref{lemma-ko-3-4}.

If $\{A',A''\}$ satisfies (\ref{lemma-ko-3-4-1}) or (\ref{lemma-ko-3-4-2}) of Lemma \ref{lemma-ko-3-4}, then we can choose a genus two handlebody component $\bar{V}_1$ in $V_1- (\cup T_i)$ which meets only two annuli of $V_1\cap (\cup T_i)$ in its closure. Moreover, the other components of $V_1-(\cup T_i)$ are all solid tori. So $G$ preserves $\bar{V}_1$, and $G$ also preserves the big solid torus $V_1-\bar{V}_1$ if $\{A',A''\}$ satisfies (\ref{lemma-ko-3-4-1}) of Lemma \ref{lemma-ko-3-4} or the set of two big solid tori $V_1-\bar{V}_1$ if $\{A',A''\}$ satisfies (\ref{lemma-ko-3-4-2}) of Lemma \ref{lemma-ko-3-4}. So we get $G\cong \mathbb{Z}_2$ or $\mathbb{D}_2$ using Lemma \ref{lemma-gm-an-1} or Lemma \ref{lemma-gm-an-2}.

If $\{A',A''\}$ satisfies (\ref{lemma-ko-3-4-3}) of Lemma \ref{lemma-ko-3-4}, then $V_1- (\cup T_i)$ consists of a genus two handlebody and three or more solid tori (one of the solid tori, say $\tilde{V}_1$ is directly induced by (iii) of Figure \ref {fig-lemma-ko-3-4}, and the others are induced from the orbit of the solid torus between $A_1$ and $A_2$ by $G$.) Since $\tilde{V}_1$ is the only solid torus component of $V_1- (\cup T_i)$ which meets only one annulus of $V_1\cap (\cup T_i)$ in its closure, say $A$, and each of the other solid tori meets two annuli of $V_1\cap (\cup T_i)$ in its closure, $G$ preserves $\tilde{V}_1$. Since $V'=V_1-\tilde{V}_1$ is a genus two handlebody and preserved by $G$, we get $G\cong \mathbb{Z}_2$ or $\mathbb{D}_2$ by Lemma \ref{lemma-gm-an-1} by the decomposition $V_1=V'\cup_{A}\overline{\tilde{V}_1}$.

\Case{3} The parallel class of $V_1\cap (\cup T_i)$ consists of three elements. Then $V_1\cap (\cup T_i)$ cuts $V_1$ into a genus two handlebody $\bar{V}_1^3$ and solid tori (see Figure \ref{fig-threeannuli}.) In particular, there is only one solid torus component of $V_1- (\cup T_i)$ which meets only one annulus of $V_1\cap (\cup T_i)$ as $V_1^1$ of Figure \ref{fig-threeannuli}.  Let this solid torus be $V_1^1$, then $V_1^1$ is preserved by $G$. Since the only genus two handlebody component, say $V_1^3$, is also preserved by $G$, $G$ preserves the big solid torus $V_1-(V_1^1\cup V_1^3)$. So we get $G\cong \mathbb{Z}_2$ or $\mathbb{D}_2$ using the same arguments in Lemma \ref{lemma-3moreannulicp}.
\end{proof}

\bibliographystyle{plain}

\end{document}